*А.В. Гасников*[1]
[1] Московский физико-технический институт (государственный университет)

## МЕТОД НЬЮТОНА – КАНТОРОВИЧА - ВЬЕТОРИСА

Пусть в банаховом пространстве $X$ действуют всюду определенные операторы $A(x)$ и $B_{n-1}(x)$, $n \in \mathbb{N}$, значения которых также принадлежат $X$. Дано уравнение

$$x = A(x). \qquad (1)$$

Рассмотрим итерационный процесс решения уравнения (1)

$$x_n = B_{n-1}(x_n),\ n \in \mathbb{N}. \qquad (2)$$

Частными случаями итерационного процесса (2) являются: 1) метод сжимающих отображений Каччопполи-Банаха $B_{n-1}(x) \equiv A(x_{n-1})$; 2) модифицированный метод Ньютона-Канторовича $B_{n-1}(x) = A'(x_0)x - A'(x_0)x_{n-1} + A(x_{n-1})$; 3) метод Ньютона-Канторовича $B_{n-1}(x) = A'(x_{n-1})x - A'(x_{n-1})x_{n-1} + A(x_{n-1})$. Исходя из результатов Л. Вьеториса 30-ых годов XX века, в работе [1] было предложено следующие обобщение методов 1) – 3):

**4)** $\|A(x_{n-1}) - B_{n-1}(x_{n-1})\| = 0$; **5)** 4) и $\|A'(x_0) - B'_{n-1}(x_0)\| = 0$; **6)** 4) и $\|A'(x_{n-1}) - B'_{n-1}(x_{n-1})\| = 0$.

Цель работы – обобщить результаты [1] на случай, когда 4) – 6) выполняются лишь приближенно. Необходимость такого обобщения связана, во-первых, с практической реализацией [2] (ошибки округления), во-вторых, в случае некорректности операторного уравнения (2), с регуляризацией [2] и, в-третьих, с приложениями к методу малого параметра [3]. Заметим, что для 1) – 3) подобного рода обобщение имеется в [2 - 4].

**Обозначения.** $r_n = \|x_{n+1} - x_n\|$, $R_n = \sum_{k=0}^{n} r_k$, $R = \lim_{n\to\infty} R_n$, $X := \bar{B}_R(x_0)$, $\tilde{r}_n = \|x_n - x_0\|$.

**Утверждение 1.** *Пусть 1)* $R < +\infty$; *2) если* $x^* = \lim_{n\to\infty} x_n$, *то* $A(x)$ *- непрерывен в точке* $x^*$; *3)* $\lim_{n\to\infty} \|A(x_n) - B_{n-1}(x_n)\| = 0$, *тогда* $x^* = A(x^*)$ *и* $\|x_n - x^*\| \le R - R_{n-1}$.

**Предположения.** *1)* $\exists\, x_0 \in X$: $\|A(x_0) - x_0\| \le \varepsilon$ *и* $\|A(x_0) - B_0(x_0)\| \le \varepsilon_0$.

*2)* $\exists\, M : \forall\, x, y \in X \to \|A(x) - A(y)\| \le M\|x - y\|$ *и*

$\forall\, n \in \mathbb{N}\ \exists\, M_{n-1} : \forall\, x, y \in X \to \|B_{n-1}(x) - B_{n-1}(y)\| \le M_{n-1}\|x - y\|$.

*3)* $\forall\, n \in \mathbb{N}\ \exists\, \varepsilon_n$ : $\|A(x_n) - B_n(x_n)\| \le \varepsilon_n$ *(метод сжимающих отображений).*

*4)* $\forall\, x \in X, n \in \mathbb{N}\ \exists\, A'(x), B'_{n-1}(x)$ *- производные по Гато.*

*5)* $\exists\, K : \forall\, x, y \in X \to \|A'(x) - A'(y)\| \le K\|x - y\|$ *и*

$\forall\, n \in \mathbb{N}\ \exists\, K_{n-1} : \forall\, x, y \in X \to \|B'_{n-1}(x) - B'_{n-1}(y)\| \le K_{n-1}\|x - y\|$.

*6)* $\forall\, n \in \mathbb{N}\ \exists\, \sigma_{n-1}$ : $\|A'(x_{n-1}) - B'_{n-1}(x_{n-1})\| \le \sigma_{n-1}$ *(метод Ньютона).*

*7)* $\forall\, n \in \mathbb{N}\ \exists\, \gamma_{n-1}$ : $\|A'(x_0) - B'_{n-1}(x_0)\| \le \gamma_{n-1}$ *(модифицированный метод Ньютона).*

**Теорема.** *Пусть справедливы предположения*

*1) 1 - 2, тогда* $r_0 \le M_0 r_0 + \varepsilon + \varepsilon_0$.

*2) 1 - 3, тогда* $r_n \le M_n r_n + (M + M_{n-1}) r_{n-1} + \varepsilon_{n-1} + \varepsilon_n$, $n \in \mathbb{N}$.

*3) 1 - 6, тогда* $r_n \le M_n r_n + 1/2\left(K + K_{n-1}\right) r_{n-1}^2 + \sigma_{n-1} r_{n-1} + \varepsilon_{n-1} + \varepsilon_n$, $n \in \mathbb{N}$.

*4) 1 - 5, 7 тогда* $\tilde{r}_n \le M_{n-1}\tilde{r}_n + 1/2\left(K + K_{n-1}\right)\tilde{r}_{n-1}^2 + \gamma_{n-1}\tilde{r}_{n-1} + \varepsilon + \varepsilon_{n-1}$, $n \in \mathbb{N}$ *и*

$$r_n \le M_n r_n + 1/2\left(K + K_{n-1}\right) r_{n-1}^2 + \left(\gamma_{n-1} + \left(K + K_{n-1}\right)\tilde{r}_{n-1}\right) r_{n-1} + \varepsilon_{n-1} + \varepsilon_n,\ n \in \mathbb{N}.$$

**Предположения.** 8) $\sum_{n=1}^{\infty} \varepsilon_n < +\infty$, $\exists M_* : \forall n \in \mathbb{N}\ \exists M_{n-1} \le M_*$.

9) $M_* < 1$, $q = \left(M + M_*\right)/\left(1 - M_*\right) < 1$.

10) $\exists K_* : \forall n \in \mathbb{N}\ \exists K_{n-1} \le K_*$.

Из теоремы и утверждения 1 следует, что при предположениях 1 – 3, 8, 9 уравнение (1) имеет решение, к которому сходится итерационный процесс (2). Если, дополнительно, справедливы предположения 4 - 6, 10 или 4, 5, 7, 10, то неравенства (в худшем случае - равенства) п.2 – 4 теоремы могут быть единообразно записаны в виде (3).

**Утверждение 2.** *Пусть* $r_0 \ge 0$; $\eta \ge 0$; $\forall n \in \mathbb{N} \to 0 \le \lambda_{n-1} \le \lambda < 1$, $\rho_{n-1} \ge 0$,

$$r_n = \eta r_{n-1}^2 + \lambda_{n-1} r_{n-1} + \rho_{n-1}. \tag{3}$$

*1) Если* $\forall n \in \mathbb{N} \to 4\eta\rho_{n-1} < \left(1 - \lambda_{n-1}\right)^2$ *и* $\exists C \ge r_0$:

$\sup\limits_{n \in N}\left(1 - \lambda_{n-1} - \sqrt{\left(1 - \lambda_{n-1}\right)^2 - 4\eta\rho_{n-1}}\right)\Big/2\eta \le C \le \inf\limits_{n \in N}\left(1 - \lambda_{n-1} + \sqrt{\left(1 - \lambda_{n-1}\right)^2 - 4\eta\rho_{n-1}}\right)\Big/2\eta$, *то*

$\forall n \in \mathbb{N} \to r_{n-1} \le C$ *(это неравенство даёт оценку сверху* $\tilde{r}_n$ *в п.4 теоремы).*

*2) Если* $\exists C_1 \in [0,1]$, $C_2 \in \left[0, \left(1 - C_1\right)^2/4\eta\right]$, $\tilde{\rho}_0 \in \left[\eta r_0^2 + \lambda_0 r_0 + \rho_0/C_\rho\right]$, *где*

$C_\rho = \left(1 - C_1 + \sqrt{\left(1 - C_1\right)^2 - 4\eta C_2}\right)\Big/2\eta C_2$: $\forall n \in \mathbb{N} \to \lambda_{n+1} \le C_1\, \rho_{n+1}/\rho_n$, $\rho_n \le C_2\, \rho_{n+1}/\rho_n$, *то*

$$\forall n \in \mathbb{N} \to \rho_n \le r_n \le C_\rho \rho_n. \tag{4}$$

*3) Если* $\exists\ \chi \in [0,1]$, $\mu \in [0, 1/\lambda - 1]$, $\tilde{\lambda}_0$, $C_\mu$: $(1 + \mu)\tilde{\lambda}_0 C_\mu \ge \eta r_0^2 + \lambda_0 r_0 + \rho_0$,

$\eta\tilde{\lambda}_0 C_\mu \le (1 - \chi)\mu\lambda_1$ *и* $\forall n \in \mathbb{N} \to \rho_n \le \chi\mu C_\mu \tilde{\lambda}_0 \prod_{k=1}^{n}\lambda_k$, $\eta C_\mu \tilde{\lambda}_0 \prod_{k=1}^{n}\lambda_k \le (1 - \chi)\mu\lambda_{n+1}$, *то*

$$\forall n \in \mathbb{N} \to r_0 \prod_{k=0}^{n-1}\lambda_k \le r_n \le C_\mu (1 + \mu)^n\, \tilde{\lambda}_0 \prod_{k=0}^{n-1}\lambda_k. \tag{5}$$

*4) Если* $\eta r_0 < 1$, $\exists\ \mu \ge 0, \chi \in [0,1]$: $\forall n \in \mathbb{N} \to \lambda_{n-1} \le \chi\mu\left(\eta r_0\right)^{2^{n-1}}$, $\eta\rho_{n-1} \le (1 - \chi)\mu\left(\eta r_0\right)^{2^n}$, *то*

$$\forall n \in \mathbb{N} \to \left(\eta r_0\right)^{2^{n-1}}\big/\eta \le r_{n-1} \le (1 + \mu)^{n-1}\left(\eta r_0\right)^{2^{n-1}}\big/\eta. \tag{6}$$

**Замечание 1.** Пусть $\forall n \in \mathbb{N} \to 4\eta\rho_{n-1} < \left(1 - \lambda_{n-1}\right)^2$ и $\left\{\lambda_n\right\}_{n=0}^{\infty}$, $\left\{\rho_n\right\}_{n=0}^{\infty}$ - невозрастающие последовательности, тогда $\forall n \in \mathbb{N} \to r_{n-1} \le \max\left\{r_0, \left(1 - \lambda_0 + \sqrt{\left(1 - \lambda_0\right)^2 - 4\eta\rho_0}\right)\Big/2\eta\right\}$.

Утверждение 2 вкупе с теоремой решают вопрос о скорости сходимости итерационного процесса (2) (см. формулы (4) – (6)) при различных предположениях относительно $B_{n-1}(x)$, $n \in \mathbb{N}$.

**Замечание 2.** Пусть $X' = C^k(T) \subset X'' = C(T) \subset X^T$, где $T$ - компакт в $\mathbb{R}^n$, $L : X' \to X''$ - линейный дифференциальный оператор на $X'$, $R : X' \to X^\gamma$ - линейный оператор граничных условий, $r_n(t) = \left\| x_{n+1}(t) - x_n(t) \right\|$. Рассмотрим дифференциальное уравнение

$$Lx(\cdot)(t) = A\left(x(t)\right) \text{ при однородных граничных условиях } Rx(\cdot) = 0. \tag{7}$$

При определенных условиях (см. [5]) (7) можно переписать как

$$x(t) = L_R^{-1} A\big(x(\cdot)\big) = \int_T G(t,s) A\big(x(s)\big)\, ds \,. \tag{8}$$

Как правило, $L_R^{-1}$ является вполне непрерывным интегральным оператором, которому соответствует операторная функция Грина $G(t,s)$ [5].

Для уравнения (8) формулировка теоремы останется такой же, только неравенства следует переписать. Продемонстрируем это на примере п.3 теоремы

$$r_n(t) \le \int_T \|G(t,s)\| \big( M_n r_n(s) + 1/2\left(K + K_{n-1}\right) r_{n-1}^2(s) + \sigma_{n-1} r_{n-1}(s) + \varepsilon_{n-1} + \varepsilon_n \big)\, ds \,.$$

Заметим, что для сходимости итерационного процесса (2) для уравнения (8) предположение 9, не нужно. Зато желательно, чтобы оператор $\int_T \|G(t,s)\| z(s) ds$ - был оператором типа Вольтера. Сказанное выше легко представить при $k=1$, $n=1$, $T=\left[0,T\right]$, $L = d/dt$.

**Замечание 3.** Пусть $B_1$, $B_2$ - банаховы пространства, $P: B_1 \to B_2$, $\Gamma: B_1 \times B_2 \to B_1$, $\Gamma(x,0) \equiv 0$. Рассмотрим уравнение $P(x)=0$, которое можно переписать как

$$x = x - \Gamma\big(x, P(x)\big)\,. \tag{9}$$

Уравнение (9) имеет вид (1), поэтому можно использовать уже полученные результаты.

Работа поддержана грантами РФФИ № 01-01-00942, офи-05-01-08045 и программой государственной поддержки ведущих научных школ НШ-1843.2003.1.

**Литература**